\documentclass{amsart}

\usepackage{amssymb}

\numberwithin{equation}{section}

\newtheorem{thm}[equation]{Theorem}
\newtheorem{cor}[equation]{Corollary}
\newtheorem{lem}[equation]{Lemma}
\newtheorem{prop}[equation]{Proposition}
\newtheorem{rem}[equation]{Remark}

\newtheorem{defn}[equation]{Definition}

\newcommand{\comment}[1]{}
\newcommand{\Ko}{K^\circ}
\newcommand{\Koo}{K^{\circ\circ}}

\renewcommand{\epsilon}{\varepsilon}
\renewcommand{\phi}{\varphi}

\let\cal\mathcal

\def\11{{\mathbf 1}}

\def\NN{{\mathbf N}}

\def\cA{{\mathcal A}}

\def\cL{{\mathcal L}}

\def\cO{{\mathcal O}}

\def\cR{{\mathcal R}}

 \def\cA{{\cal A}}
 
 \def\cR{{\cal R}}

  \def\bC{\mathbb C}
 
 \def\bQ{\mathbb Q}
 
 \def\bR{\mathbb R}
 \def\bN{\mathbb N}

\begin{document}

\title[Real closed fields with analytic structure]{Real closed fields with nonstandard and standard analytic structure}
\author[Cluckers]{R.~Cluckers$^\mathrm{1}$}
\address{Katholieke Universiteit Leuven, Departement wiskunde,
Celestijnenlaan 200B, B-3001 Leu\-ven, Bel\-gium. Current address:
\'Ecole Normale Sup\'erieure, D\'epartement de
ma\-th\'e\-ma\-ti\-ques et applications, 45 rue d'Ulm, 75230 Paris
Cedex 05, France} \email{cluckers@ens.fr}
\urladdr{www.wis.kuleuven.ac.be/algebra/Raf/}

\thanks{$^\mathrm{1}$ The author
has been supported as a postdoctoral fellow by the Fund for
Scientific Research - Flanders (Belgium) (F.W.O.) and by the
European Commission - Marie Curie European Individual Fellowship
with contract number HPMF CT 2005-007121 during the preparation of
this paper.}

\thanks{$^\mathrm{2}$ The authors have been
supported in part by NSF grant
DMS-0401175. They also thank the University of Leuven and the Newton Institute (Cambridge) for their support and hospitality.}

\author[Lipshitz]{L.~Lipshitz$^\mathrm{2}$}
\address{Department of Mathematics, Purdue University, 150 North University Street, West Lafayette IN 47907-2067, USA}
\email{lipshitz@math.purdue.edu}
\urladdr{www.math.purdue.edu/$\thicksim$lipshitz/}

\author[Robinson]{Z.~Robinson$^\mathrm{2}$}
\address{Department of Mathematics, East Carolina University, Greenville NC 27858-4353, USA}
\email{robinsonz@mail.ecu.edu}





\begin{abstract}

We consider the ordered field which is the completion of the Puiseux
series field over $\bR$ equipped with a ring of analytic functions
on $[-1,1]^n$ which contains the standard subanalytic functions as
well as functions given by $t$-adically convergent power series,
thus combining the analytic structures from \cite{DD} and
\cite{LR3}. We prove quantifier elimination and $o$--minimality in
the corresponding language. We extend these constructions and
results to rank $n$ ordered fields $\bR_n$ (the maximal completions
of iterated Puiseux series fields). We generalize the example of
Hrushovski and Peterzil \cite{HP} of a sentence which is not true in
any $o$--minimal expansion of $\bR$ (shown in \cite{LR3} to be true
in an $o$--minimal expansion of the Puiseux series field)  to a
tower of examples of sentences $\sigma_n$, true in $\bR_n$, but not
true in any $o$--minimal expansion of any of the fields
$\bR,\bR_1,\ldots,\bR_{n-1}$.

\end{abstract}

\date{}
\maketitle

\section{Introduction}
In \cite{LR3} it is shown that the ordered field $K_1$ of Puiseux
series in the variable $t$ over $\bR$, equipped with a class of
$t$--adically overconvergent functions such as
$\sum_n{(n+1)!(tx)^n}$ has quantifier elimination and is
$o$--minimal in the language of ordered fields enriched with
function symbols for these functions on $[-1,1]^n$.  This was
motivated (indirectly) by the observation of Hrushovski and
Peterzil, \cite{HP}, that there are sentences true in this structure
that are not satisfiable in any $o$--minimal expansion of $\bR$.
This in turn was motivated by a question of van den Dries.  See
\cite{HP} for details.

In \cite{DMM1} it was observed that if $K$ is a maximally complete,
non--archimedean real closed field with divisible value group, and
if $f$ an element of $\bR[[\xi]]$ with radius of convergence $>1$,
then $f$ extends naturally to an ``analytic'' function $I^n\to K$,
where $I=\{x\in K\colon -1\leq x\leq 1\}$. Hence if $\cA$ is the
ring of real power series with radius of convergence $>1$ then $K$
has $\cA$--analytic structure i.e.~this extension preserves all the
algebraic properties of the ring $\cA$. In particular the real
quantifier elimination of \cite{DD} works in this context so $K$ has
quantifier elimination, is $o$--minimal in the analytic field
language, and is even elementarily equivalent to $\bR$ with the
subanalytic structure. See \cite{DMM2} and \cite{DMM3} for
extensions.

In Sections \ref{Notation&QE} and \ref{omin} below we extend the
results of \cite{LR3} by proving quantifier elimination and
$o$--minimality for $K_1$ in a larger language that contains the
overconvergent functions together with the usual analytic functions
on $[-1,1]^n$.  In Section \ref{extensions} we extend these results
to a larger class of non-archimedean real-closed fields, including
fields $\bR_m$, of rank $m=1,2,3,\ldots$ , and in Section \ref{HP}
we show that the idea of the example of \cite{HP} can be iterated so
that for each $m$ there is a sentence true in $\bR_m$ but not
satisfiable in any $o$--minimal expansion of $\bR = \bR_0, \bR_1,
\cdots, \bR_{m-1}$.  In a subsequent paper we will give a more
comprehensive treatment of both Henselian fields with analytic
structure and real closed fields with analytic structure,
see~\cite{CL}.

\section{Notation and Quantifier Elimination}\label{Notation&QE}
In this section we establish notation and prove quantifier elimination (Theorem \ref{QE}) for the field of Puiseux series over $\bR$ in a language which contains function symbols for all the standard analytic functions on $[-1,1]^n$ and all the $t$--adically overconvergent functions on this set.
\begin{defn}
\begin{eqnarray*}
K_1&: =&\bigcup_n\bR ((t^{1/n})), \text{ the field of Puiseux series over $\bR$,}\\
K&: =& \widehat{K}_1, \text{ the $t$--adic completion of }K_1.\
\end{eqnarray*}
$K$ is a real closed, nonarchimedean normed field.  We shall use $\|\cdot\|$ to denote the (nonarchimedean) $t$--adic norm on $K$, and $<$ to denote the order on $K$ that comes from the real closedness of $K$.  We will use $|\cdot|$ to denote the corresponding absolute value, $|x| = \sqrt{x^2}$.
\begin{eqnarray*}
K^\circ&:=&\{x\in K\colon \|x\|\leq 1\},\text{ the finite elements of }K\\
K^{\circ\circ}&:=&\{x\in K\colon \|x\| < 1\},\text{ the infinitesimal elements of }K\\
K_{alg}&:=&K[\sqrt{-1}],\text{ the algebraic closure of }K\\
\xi &=&(\xi_1,\ldots,\xi_n)\\
A_{n,\alpha}&:=&\{f\in\bR[[\xi]]\colon\text{ radius of convergence of
}f>\alpha\}, 0< \alpha \in \bR\\
{\cR}_{n,\alpha}&: =& A_{n,\alpha}\widehat\otimes_{\bR} K =
(A_{n,\alpha}\otimes_{\bR} K)\;\widehat{}\, , \mbox{ where
$\,\widehat{}\,$ stands for the $t$-adic completion}
\\
{\cR^\circ}_{n,\alpha}&: =& A_{n,\alpha}\widehat\otimes_{\bR} \Ko =
(A_{n,\alpha}\otimes_{\bR} \Ko)\;\widehat{}
\\
{\cR}_n&:=&\bigcup_{\alpha>1}{\cR}_{n,\alpha}\\
\cR&: =&\bigcup_n\bigcup_{\alpha>1}{\cR}_{n,\alpha}\\
I&: =&[-1,1] = \{x\in{K}\colon |x| \leq{1}\}.
\end{eqnarray*}
\end{defn}
\begin{rem}
(i) $K$ has $\cR$--analytic structure --  in ~\cite{DMM1} it is
explained how the functions of $A_{n,\alpha}$ are defined on $I^n$.
For example, if $f \in A_{1,\alpha}$, $\alpha > 1$, $a \in
\bR\cap[-1,1]$, $\beta \in \Koo$, then $f(a+\beta) :=
\sum_n{f^{(n)}(a){\beta^n \over n!}}$.  The extension to functions
in ${\cR}_{n,\alpha}$ is clear from the completeness of $K$.  We
define these functions to be zero outside $[-1,1]^n$. This extension
also naturally works for maximally complete fields and fields of
LE-series,
see also \cite{DMM2} and \cite{DMM3}. 

(ii) $\cR$ contains all the ``standard'' real analytic functions on $[-1,1]^n$ and
all the $t$--adically overconvergent functions in the sense of \cite{LR3}.

(iii) The elements of $A_{n,\alpha}$ in fact define complex analytic functions on the complex polydisc
$\{x\in \bC \colon |x_i| \leq \alpha \text{ for } i=1,\dots,n\}$, and hence the elements of ${\cR}_{n,\alpha}$ define ``$K_{alg}$--analytic" functions on the corresponding $K_{alg}$--polydisc $\{x\in K_{alg} \colon |x_i| \leq \alpha \text{ for } i=1,\dots,n\}$.

(iv) We could as well work with $K_1$ instead of $K$.  Then we must replace ${\cR}_{n,\alpha}$ by
$$
{\cR'}_{n,\alpha}:= \bigcup_m A_{n,\alpha}\widehat\otimes_{\bR}  \bR((t^{1/m})) =   \bigcup_m \big(
(A_{n,\alpha}\otimes_{\bR}  \bR((t^{1/m})))\; \widehat \;\:\big).
$$

(v) If $\beta \leq \alpha$ and $m \leq n$ then $\cR_{\alpha, m} \subset \cR_{\beta,n}$.
\end{rem}

\begin{lem}
Every nonzero $f\in{\cR}_{n,\alpha}$ has a unique representation
\begin{equation*}
f=
\sum_{i\in I\subset \bN}
 f_i t^{\gamma_i}
\end{equation*}
where the $f_i\in A_{n,\alpha}$, $f_i\neq 0$,
$\gamma_i\in\bQ$,
 the
$\gamma_i$ are increasing, $0\in I$, and, either $I$ is finite of
the form $\{0,\ldots,n\}$, or $\gamma_i\to\infty$ as $i\to\infty$
and $I=\NN$. The function $f$ is a unit in ${\cR}_{n,\alpha}$
exactly when $f_0$ is a unit in $A_{n,\alpha}$.
\end{lem}
\begin{proof}Observe
that if $a\in{K}, a \neq 0$ then $a$ has a unique representation
$a=\sum_{i\in I\subset \bN}a_{i}t^{\gamma_i}$, where the $0\neq a_i
\in {\bR}, \gamma_i \in {\bQ}$ and, either $I$ is finite, or
$\gamma_i \to \infty$.  If $f_0$ is a unit in $A_{n,\alpha}$, then
$f\cdot f_0^{-1}\cdot t^{-\gamma_0}  = 1 + \sum_{i=1}^\infty{f_i
\cdot f_0^{-1}\cdot t^{\gamma_i - \gamma_0}}$ and $\gamma_i -
\gamma_0 > 0$.
\end{proof}

\begin{defn}
\begin{itemize}
\item[(i)] In the notation of the previous lemma, $f_0$ is called
the \emph{top slice} of $f$.\
\item[(ii)] We call $f$ \emph{regular in $\xi_n$ of degree $s$ at $a\in [I\cap\bR]^n$} if, in the classical sense,
$f_0$ is regular in $\xi_n$ of degree $s$ at $a$.\
\item[(iii)]\ We shall abuse notation and use $\|\cdot\|$ to denote the $t$--adic norm on $K$, and the
corresponding gauss--norm on  ${\cR}_{n,\alpha}$, so, with $f$ as in the above lemma, $\|f\|=\|t^{\gamma_0}\|$.\
\end{itemize}
\end{defn}

The standard Weierstrass Preparation and Division Theorems for $A_{n,\alpha}$ extend to corresponding theorems for ${\cR}_{n,\alpha}$.

 \begin{thm}\label{WPT}(Weierstrass Preparation and Division).
If $f \in {\cR}_{n,\alpha}$ with $\|f\|=1$ is regular in $\xi_n$ of
degree $s$ at 0, then there is a $\delta \in \bR$, $\delta > 0$,
such that there are unique $A_1,\ldots,A_s,U$ satisfying
$$
f=[\xi_n^s+A_1(\xi')\xi_n^{s-1}+\ldots+A_s(\xi')]U(\xi)
$$
and
$$
A_1,\ldots,A_s,\in{\cR}_{n-1,\delta}, \text{ and } U\in{\cR}_{n,\delta} \text{ a unit}.\\
$$
Then automatically
$$
\|A_1\|,\ldots,\|A_s\|,\|U\|\leq 1,\quad
\|A_1(0)\|,\ldots,\|A_s(0)\|<1,\ \text{ and }\ \|{U}(0)\|=1.
$$

Furthermore, if $g\in{\cR}_{n,\alpha}$ then there are unique
$Q\in{\cR}_{n,\delta}$ and
$R_0(\xi'),\ldots,R_{s-1}(\xi')\in{\cR}_{n-1,\delta}$, satisfying
$$\|Q\|, \|R_i\| \leq \|g\|$$
and
$$
g=Qf+R_0(\xi')+ R_1(\xi')\xi_n+\ldots+ R_{s-1}(\xi')\xi^{s-1}_n.
$$
\end{thm}
\begin{proof}
We may assume that $f = \sum_{\gamma \in I} f_{\gamma}t^\gamma,$
where  $I\subset{\bQ}^+,0\in I,I=I+I$ and $I$ is well ordered.  (We do not require that  $f_\alpha \neq 0$ for all $\alpha \in I,$ but we do require $f_0 \neq 0.$)  We prove the Preparation Theorem.  The proof of the Division Theorem is similar.
We shall produce, inductively on $\gamma \in I,$ monic polynomials $P_\gamma[\xi_n]$ with coefficients from ${\cR}_{n-1,\delta},$ and units $U_\gamma \in {\cR}_{n,\delta}$ such that, writing $\gamma'$ for the successor of $\gamma$ in $I$, we have
$$
f \equiv P_\gamma \cdot U_\gamma \text{   mod   } t^{\gamma'}
$$
and if $\gamma < \beta$
$$
P_\gamma \equiv P_\beta \text{  and  } U_\gamma \equiv U_\beta  \text{   mod   } t^{\gamma'}.
$$
Using \cite{GR}, Theorem II.D.1 (p.80), or the proof on pp. 142-144 of \cite{ZS}, we see that there is a $0<\delta \leq \alpha$ such that for every $g \in A_{n,\delta}$ the
Weierstrass data on dividing $g$ by $f_0$ are in $A_{n,\delta}.$

$P_0$ and $U_0$ are the classical Weierstrass data for $f_0$, i.e. $f_0=U_0P_0$, where $P_0 \in A_{n-1,\delta}[\xi_n]$ is monic of degree $s$, and $U \in A_{n,\delta}$ is a unit.  Suppose $P_\gamma$ and $U_\gamma$ have been found.  Then
$$
f \equiv P_\gamma \cdot U_\gamma \text {  mod  } t^{\gamma'}
$$
so we have
$$
{U_{\gamma}}^{-1} \cdot f \equiv P_\gamma + g_{\gamma'}t^{\gamma'} + o(t^{\gamma'}),
$$
where $g_{\gamma'} \in A_{n,\delta}$ and we write $o(t^{\gamma'})$ to denote terms of order $> \gamma'.$  By classical Weierstrass division we can write
$$
g_{\gamma'} = P_0 \cdot Q_{\gamma'} + R_{\gamma'},
$$
where $ Q_{\gamma'} \in A_{n,\delta}$ and $R_{\gamma'} \in A_{n-1,\delta}[\xi_n] $ has degree $< s$ in $\xi_n.$
Let
$$
P_{\gamma'} := P_\gamma + t^{\gamma'} R_{\gamma'}.
$$
Then
\begin{eqnarray*}
{U_{\gamma}}^{-1} \cdot f & =& P_\gamma + t^{\gamma'} ( P_0  Q_{\gamma'} + R_{\gamma'}) + o(t^{\gamma'})\\
&=& (P_{\gamma'} + t^{\gamma'}P_{\gamma'}Q_{\gamma'}) + t^{\gamma'}(P_0 - P_{\gamma'})Q_{\gamma'} + o(t^{\gamma'})\\
&=& P_{\gamma'}(1 + t^{\gamma'}Q_{\gamma'}) +o(t^{\gamma'}),\
\end{eqnarray*}
since $P_0 - P_{\gamma'} = o(1)$,  i.e. it has positive order. Take
$U_{\gamma'} := U_\gamma (1 + t^{\gamma'}Q_{\gamma'}).$ The
uniqueness of the $A_i$ and $U$ follows from the same induction.
\end{proof}

\begin{rem}
We remark, for use in a subsequent paper (\cite{CL}), that the argument of the previous proof works in the more general context that $I$ is a well ordered subset of the value group $\Gamma$ of a suitably complete field, for example a maximally complete field.
\end{rem}

>From the above proof or by direct calculation we have
\begin{cor}
If $g \in \cR_1$,  $\beta \in [-1,1]$ and $g(\beta)=0$ then $\xi_1-\beta$ divides $g$ in $\cR_1$.
\end{cor}

\begin{rem}
Let $f(\xi,\eta)$ be in ${\cR}_{m+n,\alpha}$. Then there are unique
$\overline{f}_{\mu}$ in ${\cR}_{m,\alpha}$ such that
$$
f(\xi,\eta)=\sum_{\mu}\overline{f}_{\mu}(\xi)\eta^{\mu}.$$
\end{rem}



The following Lemma is used to prove Theorem \ref{SNP}.
\begin{lem}\label{lemforSNP}
Let $f(\xi,\eta)
=\sum_{\mu}\overline{f}_{\mu}(\xi)\eta^{\mu}\in{\cR}_{m+n,\alpha}$.
Then the $\overline{f}_{\mu} \in {\cR}_{m,\alpha}$ and there is an
integer  $d \in \bN$,  a constant $\beta>0$, $ \beta\in\bR$ and
$g_\mu \in {\cR}^\circ_{m+n,\beta}$ for $|\mu|<d$, such that
$$
f = \sum_{|\mu| < d} {\overline{f}_\mu(\xi) g_\mu(\xi,\eta)},
$$
in $\cR_{m+n, \beta}$.
\end{lem}
\begin{proof}

We may assume that $\|f\| = 1,$ and choose a $\nu_0$ such that
$\|\overline{f}_{\nu_0} \| = 1.$ Making an $\bR$-linear change of
variables, and shrinking $\alpha$ if necessary, we may assume that
$\overline{f}_{\nu_0}$ is regular at $0$ in $\xi_m$ of degree $s$,
say.  Write $\xi'$ for $(\xi_1,\dots,\xi_{m-1})$. By Weierstrass
Division (Theorem \ref{WPT}) there is a $\beta > 0$ and there are
$Q(\xi,\eta) \in  {\cR}_{m+n,\beta}$ and $R(\xi,\eta) =
R_0(\xi',\eta)+ \cdots + R_{s-1}(\xi',\eta)\xi_m^{s-1} \in
{\cR}_{m+n-1,\beta}[\xi_m]$ such that
$$
f(\xi,\eta) = \overline{f}_{\nu_0}(\xi) Q(\xi,\eta) + R(\xi,\eta).
$$
By induction on m, we may write
$$
R_0 = \sum_{|\mu|<d}{\overline{R}_{0\mu}(\xi')
g_\mu(\xi',\eta)},
$$
for some $d \in \bN$, some $ \beta > 0$ and $g_\mu(\xi',\eta) \in
{\cR}_{m+n-1,\beta}^\circ.$ Writing $R = \sum_\nu{\overline
R_\nu(\xi)\eta^\nu}$, observe that each $\overline{R}_\nu$ is an
${\cR}_{m,\beta}^\circ-$linear combination of the
$\overline{f}_{\nu
},$ since, taking the coefficient of $\eta^\nu$
on both sides of the equation $f(\xi,\eta) =
\overline{f}_{\nu_0}(\xi) Q(\xi,\eta) + R(\xi,\eta),$ we have
$$
\overline{f}_\nu = \overline{f}_{\nu_0}\overline{Q}_\nu +
\overline{R}_\nu.
$$
Consider
\begin{eqnarray*}
f -  \overline{f}_{\nu_0} Q -
\sum_{|\mu|<d}{\overline{R}_{\mu}(\xi)g_\mu(\xi',\eta)}
&=:& S_1\xi_m + S_2{\xi_m}^2 + \cdots +S_{s-1} \xi_m^{s-1}\\
&=& \xi_m[S_1 + S_2\xi_m + \cdots + S_{s-1}\xi_m^{s-2}]\\
&=:& \xi_m \cdot S,\, \text{ say,}
\end{eqnarray*}
where the $S_i \in {\cR}_{m+n-1,\beta}^\circ.$ Again, observe that
each $\overline{S}_\nu$ is an ${\cR}_{m,\beta}^\circ-$linear
combination of the $\overline{f}_{\nu'}.$  Complete the proof by
induction on $s$, working with $S$ instead of $R$.
\end{proof}

\begin{thm}\label{SNP}(Strong Noetherian Property).
Let $f(\xi,\eta)
=\sum_{\mu}\overline{f}_{\mu}(\xi)\eta^{\mu}\in{\cR}_{m+n,\alpha}$.
Then the $\overline{f}_{\mu} \in {\cR}_{m,\alpha}$ and there is an
integer  $d \in \bN$,  a constant $\beta>0$, $ \beta\in\bR$ and
units $U_\mu(\xi,\eta) \in {\cR}^\circ_{m+n,\beta}$ for $|\mu|<d$,
such that
$$
f = \sum_{\mu\in J} {\overline{f}_\mu(\xi) \eta^\mu
U_\mu(\xi,\eta)}
$$
in $\cR_{m+n, \beta}$, where $J$ is a subset of $\{0,1,\ldots,d\}^n$.
\end{thm}
\begin{proof}
It is sufficient to show that there are an integer $d$, a set
$J\subset\{0,1,\ldots,d\}^n$, and $g_\mu \in
{\cR}_{m+n,\beta}^\circ$ such that
 \begin{equation}\label{eqfg}
 f = \sum_{\mu\in J}
{\overline{f}_\mu(\xi) \eta^\mu g_\mu(\xi,\eta)},
\end{equation}
since then, rearranging the sum if necessary, we may assume that
each $g_\mu$ is of the form $1 + h_\mu$ where $h_\mu \in
(\eta){\cR}^\circ_{m+n,\beta}.$ Shrinking $\beta$ if necessary will
guarantee that the $g_\mu$ are units. But then it is in fact
sufficient to prove (\ref{eqfg}) for $f$ replaced by
$$
f_{I_i}:=\sum_{\mu\in I_i}\overline{f}_{\mu}(\xi)\eta^{\mu}
$$
for each $I_i$ in a finite partition $\{I_i\}$ of $\bN^n$.

By Lemma \ref{lemforSNP} there is an integer  $d \in \bN$,  a
constant $\beta>0$, $ \beta\in\bR$ and $g_\mu \in
{\cR}^\circ_{m+n,\beta}$ for $|\mu|\leq d$, such that
$$
f = \sum_{|\mu| \leq d} {\overline{f}_\mu(\xi) g_\mu(\xi,\eta)}.
$$
Rearranging, we may assume for $\nu,\mu\in\{1,\ldots,d\}^n$ that
$(\bar g_{\mu})_\nu$ equals $1$ if $\mu=\nu$ and that it equals $0$
otherwise.


Focus on $f_{I_1}(\xi,\eta)$, defined as above by
$$
f_{I_1}(\xi,\eta) = \sum_{\mu\in
I_1}\overline{f}_{\mu}(\xi)\eta^{\mu}
$$
with
$$
I_1:=\{0,\ldots,d\}^n\cup\{\mu \colon \mu_i\geq d\mbox{ for all $i$}\}
$$
and note that
\begin{equation}\label{eqhatf}
f_{I_1}(\xi,\eta) = \sum_{|\mu| \leq d} {\overline{f}_\mu(\xi)
g_{\mu,I_1}(\xi,\eta)}
\end{equation}
with $g_{\mu, I_1}(\xi,\eta)\in \cR^\circ_{m+n,\beta}$ defined by
the corresponding sum
$$
g_{\mu,I_1}(\xi,\eta) = \sum_{\nu\in
I_1}\overline{g}_{\mu,\nu}(\xi)\eta^{\nu}.
$$


It is now clear that $g_{\mu,I_1}$ is of the form
$\eta^{\mu}(1+h_\mu)$ where $h_\mu\in
(\eta){\cR}^\circ_{m+n,\beta}$.


One now proceeds by noting that $f-f_{I_1}$ is a finite sum of terms
of the form $f_{I_j}$ for $j>1$ and $\{I_j\}_j$ a finite partition
of $\bN^n$  and where each $f_{I_j}$ for $j>1$ is of the form $\eta_i^\ell
q(\xi,\eta')$ where $\eta'$ is
$(\eta_1,\ldots,\eta_{i-1},\eta_{i+1},\ldots,\eta_{n})$ and $q$ is
in $\cR_{m+n-1,\beta}^\circ$. These terms can be handled by
induction on $n$.

\end{proof}


\begin{defn}
For $\gamma\in K_{alg}^\circ$ let $\gamma^\circ$ denote the closest element of $\bC$, i.e. the unique element $\gamma^\circ$ of $\bC$ such that $|\gamma - \gamma^\circ| \in K_{alg}^{\circ\circ}.$
\end{defn}

\begin{lem}\label{standpart}
Let $f \in {\cR}_1.$ If $f(\gamma)=0$ then $f_0(\gamma^\circ)=0$
($f_0$ is the top slice of $f$). Conversely, if $\beta\in\bR$ (or
$\bC$) and $f_0(\beta)=0$ there is a $\gamma\in K_{alg}$ with
$\gamma^\circ=\beta$ and $f(\gamma)=0$. Indeed, $f_0$ has a zero of
order $n$ at $\beta\in\bC$ if, and only if, $f$ has $n$ zeros
$\gamma$ (counting multiplicity) with ${\gamma}^\circ = \beta$.
\end{lem}

\begin{proof}Use Weierstrass Preparation and \cite{BGR} Proposition 3.4.1.1.
\end{proof}

\begin{cor}\label{2.8}A nonzero $f\in{\cR}_{1,\alpha}$ has only finitely many zeros in the set $\{x\in K_{alg}\colon |x|\leq\alpha\}$.
Indeed, there is a polynomial $P(x)\in K[x]$ and a unit $U(x)\in{\cR}_{1,\alpha}$ such that $f(x)=P(x)\cdot U(x)$.
\end{cor}

\begin{proof}Observe that $f_0$ has only finitely many zeros
in $\{x\in\bC\colon |x|\leq\alpha\}$, that non-real zeros occur in
complex conjugate pairs, and that $f$ is a unit exactly when $f_0$
has no zeros in this set, i.e.~when $f_0$ is a unit in
$A_{1,\alpha}$, and use Lemma \ref{standpart}.
\end{proof}

\medskip\noindent
\begin{thm}[Quantifier Elimination Theorem]\label{QE}
Denote by  ${\cL}$  the language $\langle +,\cdot,^{-1},0,1,<,{\cR}\rangle$ where the functions in ${\cR}_n$
are interpreted to be zero outside $I^n$.
Then $K$ admits quantifier elimination in ${\cL}$.
\end{thm}

\begin{proof}
This is a small modification of the real quantifier elimination of
\cite{DD} as in \cite{DMM1}, using the Weierstass Preparation
Theorem and the Strong Noetherian Property above. Crucial is that,
in Theorems \ref{WPT} and \ref{SNP}, as in \cite{DD},
$\beta$ and $\delta$ are positive real numbers so one can use the
compactness of $[-1,1]^n$ in $\bR^n$.
\end{proof}

\section{$o$--minimality}\label{omin}

In this section we prove the $o$--minimality of $K$ in the language $\cL$. Let $\alpha > 1.$ As we remarked above, each $f\in A_{n,\alpha}$ defines a function from the poly-disc $(I_{\bC, \alpha})^{n} \to\bC$, where $I_{\bC, \alpha}:=\{x\in\bC\colon |x|\leq \alpha \}$, and hence each $f\in {\cR}_{n,\alpha}$
defines a function from $(I_{K_{alg}, \alpha})^{n} \to K_{alg}$, where $I_{K_{alg}, \alpha}:=\{x \in K_{alg} \colon |x| \leq \alpha \}$.
In general $A_{n,\alpha}$ is not closed under composition.
However, if $F(\eta_1,\ldots,\eta_m)\in {\cR}_{m,\alpha}$, $G_j(\xi)\in {\cR}_{n,\beta}$ for $j=1,\ldots,m$ and $|G_j(x)|\leq \alpha$ for all $x\in (I_{K_{alg}, \beta})^n$, then $F(G_1(\xi),\ldots,G_m(\xi))\in {\cR}_{n,\beta}$.  This is clear if $F \in A_{n,\alpha}$ and the $G_j \in A_{n,\beta}.$  The general case follows easily.

For $c,r\in K,\ r>0$, we denote the ``closed interval'' with center $c$ and radius $r$ by
$$
I(c,r) := \{x\in K\colon |x-c|\leq r\}
$$
and for $c,\delta,\epsilon\in K$, $0 < \delta < \epsilon$, we denote the ``closed annulus'' with center $c$, inner radius $\delta$ and outer radius $\epsilon$ by
$$
A(c,\delta,\epsilon) :=\{x\in K\colon\delta\leq |x-c|\leq\epsilon\}.
$$
On occasion we will consider $I(c,r)$ as a disc in $K_{alg}$ and $A(c,\delta,\epsilon)$ as an annulus in $K_{alg}$, replacing $K$ by $K_{alg}$ in the definitions.
No confusion should result.
Note that these discs and annuli are defined in terms of the real-closed order on $K$, not the non--archimedean absolute value $\|\cdot\|$, and hence are not discs or annuli in the sense of \cite{BGR},\cite{LR3} or \cite{FP},  which we will refer to as \emph{affinoid discs} and \emph{affinoid annuli.}
For $I=I(c,r)$, $A=A(c,\delta,\epsilon)$ as above, we define the rings of analytic function on $I$ and $A$ as follows:
\begin{eqnarray*}
{\cO}_I&:=&\left\{f\left({x-c\over r }\right)\colon f\in {\cR}_1\right\}\\
{\cO}_A&:=&\left\{g\left({\delta\over x-c }\right)+h\left({x-c\over\epsilon}\right)\colon g,h\in{\cR}_1,\ g(0) = 0\right\}.
\end{eqnarray*}
The elements of ${\cO}_I$ (respectively, ${\cO}_A$) are analytic functions on the corresponding $K_{alg}$--disc (respectively, annulus) as well.

\begin{rem}\label{prod}
(i) Elements of ${\cO}_A$ are multiplied using the relation ${\delta \over x-c}\cdot{x-c \over \epsilon} = {\delta \over \epsilon}$ and the fact that $|{\delta \over \epsilon}| <1$.  Indeed, let $g(\xi_1) = \sum_i{a_i\xi_1^i}$, $h(\xi_2) = \sum_j{b_j\xi_2^j}$.  Then, using the relation $\xi_1\xi_2 = {\delta \over \epsilon}$, we have
\begin{eqnarray*}
g\cdot h &=& \sum_{j<i}{a_ib_j\big({\delta \over \epsilon}\big)^j\xi_1^i} + \sum_{i\leq j}{a_ib_j\big({\delta \over \epsilon}\big)^{j-i}\xi_2^{j-i}}\\
&=& f_1(\xi_1) + f_2(\xi_2).\
\end{eqnarray*}
If $g$, $h \in A_{1,\alpha}$ and ${\delta \over \epsilon} \in \bR$
then $f_1$, $f_2 \in A_{1,\alpha}$, and this extends easily to the
case $g$, $h \in {\cR}_{1,\alpha}$ and ${\delta \over \epsilon} \in
\Ko$.  Lemma \ref{3.3} will show that in fact the only case of an
annulus that we must consider is when ${\delta \over \epsilon} \in
\Koo$.

(ii) We define the gauss-norm on $\cO_I$ by $\|f\big({x-c \over r}\big)\| := \|f(\xi)\|$, and on $\cO_A$ by
$\|g\big({\delta \over x-c}\big) + h\big({x-c \over r}\big)\| := \max \{\|g(\xi_1)\|, \|h(\xi_2)\|\}$.  It is clear that the gauss-norm equals the supremum norm.

(iii) If $f \in \cO_I$ then $\|\big({x-c \over r}\big)f\| = \|f\|$.
If $f \in \cO_A$ then $\|\big({x-c \over \epsilon}\big)f\| \leq
\|f\|$ and if ${\delta \over \epsilon} \in \Koo$, (i.e. is infinitesimal) then
$\|\big({x-c \over \epsilon}\big)g\big({\delta \over x-c}\big)\| =
\|{\delta \over \epsilon}\| \cdot \|g\big({\delta \over x-c}\big)\| <
\|g\big({\delta \over x-c}\big)\|$ and
$ \|\big({x-c \over \epsilon}\big)h\big({x-c \over \epsilon}\big)\| =
 \|h\big({x-c \over \epsilon}\big)\|$.

 (iv) If $\|f\| < 1$ then $1-f$ is a unit in $\cO_A$.
 (In fact it is a \emph{strong unit}  --  a unit $u$ satisfying $\|1-u\| < 1$.)
\end{rem}

\begin{defn}
(i) We say that $f\in{\cO}_{I(0,1)}={\cR}_1$ \emph{has a zero close to $a \in I(0,1)$} if $f$, as a $K_{alg}$--function defined on the $K_{alg}$--disc $\{|x|\colon |x|\leq\alpha\}$ for some $\alpha>1$ has a $K_{alg}$--zero  $b$ with $a-b$ infinitesimal in $K_{alg}$.  We say f \emph{has a zero close to $I(0,1)$} if it has a zero close to $a$ for some $a \in I(0,1)$.
For an arbitrary interval $I=I(c,r)$ we say that $f=F({x-c\over r})\in{\cO}_I$ \emph{has a zero close to $a \in I$} if $F$ has a zero close to ${a-c \over r} \in I(0,1)$, and that $f$ \emph{has a zero close to $I(c,r)$} if $F$ has a zero close to $I(0,1)$.

(ii) For $0 < a,b \in \Ko$ we write $a\thicksim b$ if ${a \over b}, {b \over a} \in \Ko$ and we write $a << b$ if ${a \over b} \in \Koo$.

(iii) Let $X$ be an interval or an annulus, and let $f$ be defined
on a superset of $X$. We shall write $f\in{\cO}_X$ to mean that
there is a function $g\in {\cO}_X$ such that
$$
f|_X=g.
$$
\end{defn}

\begin{lem}\label{3.1}If $f\in{\cO}_{I(c,r)}$ has no zero close to $I(c,r)$, then there is a cover of $I(c,r)$ by finitely many closed intervals $I_j=I(c_j,r_j)$ such that ${1\over f}\in{\cO}_{I_j}$ for each $j$.
\end{lem}

\begin{proof}It is sufficient to consider the case $I(c,r)=I(0,1)$.
Cover $I(0,1)$ by finitely many intervals $I(c_j,r_j)$, $c_j,r_j\in{\bR}$ such that $f$ has no $K_{alg}$--zero in the $K_{alg}$--disc $I(c_j,r_j)$.
Finally use Corollary \ref{2.8}.
\end{proof}

\begin{rem}\label{subset}
The function $f(x) = 1+x^2$ has no zeros close to $I(0,1)$.  It is
not a unit in $\cO_{I(0,1)}$, but it is a unit in both $\cO_{I(-{1
\over 2}, {1 \over 2})}$ and $\cO_{I({1 \over 2}, {1 \over 2})}$.
The function $g(x)=x$ is not a unit in $\cO_{I(\delta, \epsilon)}$
for any $0 < \delta \in \Koo$ and $0 < \epsilon \in \Ko \setminus
\Koo$.  It is of course a unit in $\cO_{A(0,\delta, \epsilon)}$. The
function ${1 \over g} = {1 \over x} \in A(0,\delta,1)$ for all $\delta > 0$ but is not in
$\cO_I$ for $I = I({1+\delta \over 2},{1-\delta \over 2})$, for any
$\delta \in \Koo$.  Thus we see that if $X_1 \subset X_2$ are
annuli or intervals it does not necessarily follow that $\cO_{X_2}
\subset \cO_{X_1}$.  However the following are clear. If $I_1
\subset I_2$ are intervals, then $\cO_{I_2} \subset \cO_{I_1}$. If
$0 < \delta \in \Koo$, $0 < r \in \Ko \setminus \Koo$, $r < 1$, and $A= A(0,
\delta, 1)$, $I = I({1+r \over 2},{1-r \over 2})$ then $\cO_A
\subset \cO_I$.  If $A_1 \subset A_2$ are annuli that have the same
center, then $\cO_{A_2} \subset \cO_{A_1}$.  If $0 < \delta << c$
and $0 < r < {c \over \alpha}$ for some $1 < \alpha \in \bR$, and $I
= I(c,r) \subset A(0,\delta,1) = A$, then $\cO_A \subset \cO_I$.
(Writing $x = c - y$, $|y| \leq r$ we see that ${\delta \over x} =
{\delta \over {c-y}} ={{\delta/c} \over{1-y/c}} + {\delta \over
c}\sum{({y \over c})^k}= {\delta \over c} \sum{({r \over c}{y \over
r})^k}$).
\end{rem}

Restating Corollary \ref{2.8} we have

\begin{cor}\label{polunit}If $f\in{\cO}_{I(c,r)}$ there is a polynomial $P\in K[\xi]$ and a unit $U\in{\cO}_{I(c,r)}$ such that $f(\xi)=P(\xi)\cdot U(\xi)$.
\end{cor}

\begin{lem}\label{3.3}If $\epsilon < N\delta$ for some $N\in\bN$, then there is a covering of $A(c,\delta,\epsilon)$ by finitely many intervals $I_j$ such that for every $f\in{\cO}_{A(c,\delta,\xi)}$ and each $j$, $f\in{\cO}_{I_j}$.
\end{lem}

\begin{proof}Use Lemma \ref{3.1} or reduce directly to the case $\epsilon = 1$, $0 < \delta = r \in \bR$ and the two intervals $[-1,-r] = I(-{1+r \over 2}, {1-r \over 2})$ and $[r,1] = I({1+r \over 2}, {1-r \over 2})$.
\end{proof}

The following Lemma is key for proving $o$--minimality.

\begin{lem}\label{cover}Let $f\in A(c,\delta,\epsilon)$.
There are finitely many intervals and annuli $X_j$ that cover $A(c,\delta,\epsilon)$, polynomials $P_j$ and units $U_j\in{\cO}_{X_j}$ such that for each $j$ we have $f|_{X_j}=(P_j\cdot U_j)|_{X_j}$.
\end{lem}

\begin{proof}By the previous lemma, we may assume that $c=0,\ \epsilon=1$ and $\delta\in K^{\circ\circ}$ (i.e.
$\delta$ is infinitesimal, say $\delta = t^{\gamma}$ for some
$\gamma >0$). Let
$$
f(x)=g\big({\delta\over x}\big)+h(x)
 \quad\mathrm{with}\quad  g(\xi),\ h(\xi)\in{\cR}_1,\ g(0)=0,
$$
and
$$
g(\xi)=
\sum_{i\in I\subset \NN}
 t^{\alpha_i}\xi^{n_i} g_i(\xi),
$$
with $n_i>0$, $g_i(0) \neq 0$, $g_i \in A_{1,\alpha}$ for some $\alpha >1$.
Observe that
$$
xg\big({\delta\over x}\big)=\sum_{i\in I}
(t^{\alpha_i}\delta)\big({\delta\over x}\big)^{n_i - 1}
g_i\big({\delta\over x}\big).
$$
Hence (see Remark \ref{prod}) for suitable $n\in\bN$, absorbing the
constant terms into $h$, we have
$$
x^n f(x)=\overline g\big({\delta\over x}\big)+\overline h(x)
$$
where $\overline g(0)=0$ and $\|\overline g\| < \|\overline h\|$.  (For use in Section 4, below, note that this argument does not use that $K$ is complete or of rank $1$.)
Multiplying by a constant, we may assume that $\|\overline h\|=1$.
Let
$$
\overline g(\xi)= \sum_{i\in \overline I}
t^{\beta_i}\xi^{m_i}\overline g_i(\xi), \text{ with } m_i>0 \text{
and } \beta_i>0 \text{ for each } i,
$$
and
$$
\overline h(\xi)=\xi^{k_0} \overline h_0(\xi)+\sum\limits_{i\in
J\subset \NN\setminus{0}} t^{\gamma_i} \overline h_i(\xi), \text{
with } \overline h_0 (0)\neq0 \text{ and the } \gamma_i>0
\text{ increasing}.
$$
Since $\|\overline g\| = \|t^{\beta_0}\| < 1$ there is a $\delta'$ with
$\delta\leq\delta'\in K^{\circ\circ}$ (i.e.~$\|\delta'\|<1$) such
that
 $\|t^{\gamma_1}\| < \|(\delta')^{k_0}\|$ and $\|\overline g\|<\|(\delta')^{k_0}\|$.
Splitting off some intervals of the form
$I(\frac{-1-r}{2},\frac{-1+r}{2})=[-1,-r]$ or
$I(\frac{1+r}{2},\frac{1-r}{2})=[r,1]$ for $r>0$, $r\in\bR$ (on
which the result reduces to Corollary \ref{polunit} by Remark
\ref{subset}) and renormalizing, we may assume that $\overline h_0,$ is a unit
in $A_{1,\alpha}$, for some $\alpha > 1$.  So
\begin{eqnarray*}
x^n f(x)&=&\overline g\big({\delta\over x}\big)+\overline h(x)\\
&=&\overline h_0(x) x^{k_0}\Big[1+ \sum_{i=1}^\infty{\big({\delta'\over x}\big)^{k_0}{t^{\gamma_i} \over (\delta')^{k_0}}
{\overline h_i(x) \over \overline h_0(x)}} + {1 \over \overline h_0(x)}\big({1 \over \delta'}\big)^{k_0}\big({\delta'\over x}\big)^{k_0}g\big({\delta \over \delta'}{\delta'\over x}\big)\Big].\
\end{eqnarray*}
By our choice of $\delta'$ and Remark \ref{prod} the quantity in square brackets is a (strong) unit. Hence we have taken care of an annulus of the form $A(0,\delta',1)$ for some $\delta'$ with $|\delta| \leq |\delta'|$ and $\|\delta'\| < 1$.

Observe that the change of variables $y={\delta\over x}$ interchanges the sets $\{x\colon \|x\|=\|\delta\|$ and $\delta\leq |x|\}$ and $\{y\colon \|y\|=1$ and $|y|\leq 1\}$.
Hence, as above, there is a $\delta''\in K^{\circ\circ}$ with $\|\delta\| < \|\delta''\|$ and a covering of the annulus $\delta\leq |x|\leq\delta''$ by finitely many intervals and annuli with the required property.

It remains to treat the annulus $\delta''\leq |x|\leq\delta'$.
Using the terminology of \cite{LR3}, observe that on the much bigger \emph{affinoid} annulus $\|\delta''\|\leq \|x\|\leq \|\delta'\|$ the function $f$ is strictly convergent, indeed even overconvergent.
Hence, as in \cite{LR3} Lemma 3.6, on this affinoid annulus we can write
$$
f={P(x)\over x^\ell}\cdot U(x)
$$
where $P(x)$ is a polynomial and $U(x)$ is a strong unit (i.e.~$\|U(x)-1\|<1$.)
\end{proof}

\begin{cor}\label{3.5}If  $X$ is an interval or an annulus and
$f\in{\cO}_X,$  then the set $\{x \in X\colon f(x)\geq 0\}$ is
semialgebraic (i.e.~a finite union of (closed) intervals).
\end{cor}

\begin{proof}
This is an immediate corollary of \ref{polunit} and \ref{cover}
since units don't change sign on intervals and since an annulus has
two intervals as connected components.
\end{proof}

\begin{defn}For $c=(c_1,\ldots,c_n)$, $r=(r_1,\ldots,r_n)$ we define the poly--interval $I(c,r):=\{x\in K^n\colon |x_i-c_i|\leq r_i$, $i=1,\ldots,n\}$.  This also defines the corresponding polydisc in $(K_{alg})^n.$  The ring of analytic functions on this poly-interval (or polydisc) is
$$
{\cO}_{I(c,r)}:=\left\{f\left({x_1-c_1\over r_1 },\ldots,{x_n-c_n\over r_n}\right)\colon f\in{\cR}_n\right\}.
$$
\end{defn}

\begin{lem}\label{compact}Let $\alpha, \beta > 1, F(\eta_1,\ldots,\eta_m)\in{\cR}_{m,\alpha}$ and $G_j(\xi_1, \cdots, \xi_n)\in{\cR}_{n,\beta}$ with $\|G_j\| \leq 1$ for $j=1,\ldots,m$.
Let $X=\{x\in [-1,1]^n\colon |G_j(x)|\leq 1$ for $j=1,\ldots,m\}$.
There are (finitely many) $c_i =(c_{i1},\ldots,c_{in})\in\bR^n,\ \epsilon_i\in\bR$, $\epsilon_i>0$ with $\epsilon_i< |c_{ij}|$ if $c_{ij}\neq 0$  such that the (poly) intervals $I_i = I(c_i,\epsilon_i) =\{x\in K\colon |x_j-c_{ij}| <\epsilon_i$ for $j=1,\ldots,n\}$ cover $X$, $|G_j(x)|<\alpha$ for all $x\in I_i,j=1,\ldots,m$, and there are $H_i\in {\cO}_{I_i}$ such that
$$
F(G_1,\ldots,G_m)|_{I_i}=H_i|_{I_i}.
$$
\end{lem}

\begin{proof}Use the compactness of $[-1,1]^n\cap\bR^n$ and the following facts.  If $\|G_j\| = 1$ then $|G_j(x)-G_{j0}(x)|$ is infinitesimal for all $x \in \Ko$, where $G_{j0}$ is the top slice of $G_j$.  If $\|G_j\| < 1$ then $|G_j(x)| \in \Koo$ for all $x \in [-1,1]^n$.
\end{proof}

\begin{cor}\label{3.10}(i) Let $I$ be an interval, $F(\eta_1,\ldots,\eta_m)\in{\cR}_m$, and $G_j\in{\cO}_I$ for $j=1,\ldots,m$.
Then there are finitely many intervals $I_i$ covering $I$ and functions $H_i\in{\cO}_{I_i}$ such that for each $i$
$$
F(G_1,\ldots,G_m)|_{I_i}=H_i|_{I_i}
$$
(ii)\ Let $A$ be an annulus, $F(\eta_1,\ldots,\eta_m)\in{\cR}_m$, and $G_j\in{\cO}_A$ for $j=1,\ldots,m$.
Then there are finitely many $X_i$, each an interval or an annulus, covering $A$ and $H_i\in{\cO}_{X_i}$ such that for each $i$
$$
F(G_1,\ldots,G_m)|_{X_i}=H_i|_{X_i}.
$$
\end{cor}

\begin{proof}Part (i) reduces to Lemma \ref{compact} once we see that if $\|G_j\| > 1$ we can use Corollary \ref{2.8} to restrict to the subintervals of $I$ around the zeros of $G_j$ on which $|G_j| \leq C$ for some $1 < C \in \bR$. On the rest of $I$, $F(G_1,\ldots,G_m)$ is zero.

For (ii), we may assume that $A=A(0,\delta,1)$ with $\delta$
infinitesimal and that $G_j(x)=G_{j_1}({\delta\over x})+G_{j_2}(x)$.
As in (i), we may reduce to the case that $\|G_j\| \leq 1$, using
Lemma \ref{cover} instead of Corollary \ref{2.8}, and using Lemma
\ref{3.3} and Remark \ref{subset}. Apply Lemma \ref{compact} to the
functions $F$ and $G'_j(\xi_1,\xi_2)=G_{j_1}(\xi_1)+G_{j_2}(\xi_2)$.
The case $c=(0,0)$ gives us the annulus $|{\delta\over
x}|\leq\epsilon$, $|x|\leq\epsilon$ i.e.~${\delta\over\epsilon}\leq
|x|\leq\epsilon$. The case $c=(0,c_2)$ with $c_2\neq 0$ gives us
${\delta\over |x|}\leq\epsilon$ and $|x-c_2|\leq\epsilon$ with
$\epsilon<|c_2|$. This is equivalent to $|x-c_2|\leq\epsilon$ since
$\epsilon,c_1\in\bR$, and hence equivalent to $-\epsilon\leq
x-c_2\leq\epsilon$ or $c_2-\epsilon\leq x\leq c_2+\epsilon$ which is
an interval bounded away from 0. The case $c=(c_1,0)$ gives us
$|{\delta\over x}-c_1|\leq\epsilon$, $|x|\leq\epsilon$ (since
$\epsilon<|c_1|)$ which is equivalent to $|{\delta\over
x}-c_1|\leq\epsilon$ or $c_1-\epsilon\leq {\delta\over x}\leq
c_1+\epsilon$ or (considering the case $c_1>0$, the case $c_1<0$ is
similar) ${\delta\over c_1+\epsilon}\leq x\leq{\delta\over
c_1-\epsilon}$ which is part of an annulus that can be reduced to
intervals using Lemma \ref{3.3}. The case $c=(c_1,c_2)$ with both
$c_1,c_2\neq 0$ is vacuous since either $x$ or ${\delta\over x}$ is
infinitesimal on $A$ and $\epsilon < |c_1|, |c_2|$.
\end{proof}

\begin{lem}\label{3.11}Let $X$ be an interval or an annulus and let $f,g\in{\cO}_X$.
There are finitely many subintervals and subannuli $X_i \subset X$, $i=1,\ldots,\ell$ such that
$$
\{x\in X\colon |f(x)|\leq |g(x)|\}\subset \bigcup X_i\subset X,
$$
and, except at finitely many points,
$$
{f\over g}|_{X_i}\in{\cO}_{X_i}.
$$
\end{lem}

\begin{proof}We consider the case that $X$ is an interval, $I$, and may take $I= I(0,1) =[-1,1]$. We may assume by Corollary \ref{polunit} that $f$ and $g$ have no common zero.
If $g$ has no zeros close to $I(0,1)$, we are done by Lemma \ref{3.1}.
Let $\alpha_1,\ldots,\alpha_n$ be the distinct elements $\alpha$ of $[-1,1]\cap\bR$ such that $g$ has at least one zero close to (i.e.~within an infinitesimal of) $\alpha$.
Breaking into subintervals and making changes of variables we may assume that $n=1$ and that $\alpha_1=0$.
Again making a change of variables (over $K$) we may assume that $g$ has at least one zero with zero ``real" part, i.e. of the form $a = \sqrt{-1}\alpha$ for some $\alpha \in K$.
Let $N$ denote the number of zeros of $f\cdot g$ close to 0 in $I(0,1)$.
Let
$\delta = 3\cdot\max\{|x|\colon x \in K_{alg} \text{ close to }0 \text{ and } f(x)\cdot g(x)=0\}$.
If $\delta = 0$, then $a = \alpha = 0$ and there is no other zero of $f\cdot g$ close to zero in $I(0,1)$.  Then there is a $\delta' > 0$ such that for $|x| < \delta'$ we have $|g(x)| > |f(x)|$.  Then the interval $I(0,\delta')$ drops away, and on the annulus $A(0,\delta',1)$ the function $g$ is a unit.  If $\delta>0$,
we consider the interval $I(0,\delta)$ and the annulus $A(0,\delta,1)$ separately, and proceed by induction on $N$.
So suppose $\delta > 0$ and $f\cdot g$ has $N$ zeros close to 0 in $I(0,1)$.
Let $\alpha$ be as above.  If $\alpha \thicksim \delta$ then the zero $a = \sqrt{-1}\alpha$ is not close to $I(0,\delta)$ and by restricting to $I(0,\delta)$ we have reduced $N$.  If $|\alpha| << \delta$ then this zero is close to $0$ in $I(0,\delta)$, but the largest zeros (those of size $\delta/3$) are not close to $0$ in $I(0,\delta)$, and hence restricting to $I(0,\delta)$ again reduces $N$.

It remains to consider the case of the annulus  $A(0,\delta,1) = \{x\colon \delta\leq |x|\leq 1\}$ where all the zeros of $g$ are within $\delta/3$ of 0.
By Lemma \ref{cover} we may assume that $g(x)=P(x)\cdot U(x)$ where $U$ is a unit and all the zeros of $P$ are within $\delta/3$ of 0.
Let $\alpha_i$, $i=1,\ldots,\ell$ be these zeros.
For $|x|\geq\delta$ we may write ${1\over x-\alpha_i}={1\over x}{1\over 1-\alpha_i/x}={1\over x}\sum^\infty_{j=0} ({\alpha_i\over x})^j={1\over x} \sum^\infty_{j=0} ({\alpha_i\over\delta})^j ({\delta\over x})^j$ and $|{\alpha_i\over\delta}| \leq {1\over 3}$.
Hence ${f\over g}\in {\cO}_{A(0,\delta,1)}$.
This completes the case that $X$ is an interval.

The case that $X$ is an annulus is similar  --  one can cover $X$ with finitely many subannuli $X_i$ and subintervals $Y_j$ so that for each $i$ \, $g|_{X_i}$ is a unit in $\cO_{X_i}$ and for each $j$ \, $f|_{Y_j},
g|_{Y_j} \in \cO_{Y_j}$.
\end{proof}

>From Corollary \ref{3.10} and Lemma \ref{3.11} we now have by induction on terms:

\begin{prop}Let $f_1,\ldots,f_\ell$ be ${\cL}$--terms in one variable, $x$.
There is a covering of $[-1,1]$ by finitely many intervals and annuli $X_i$ such that  except for finitely many values of $x$, we have for each $i$ and $j$ that $f_i|_{X_j}\in{\cO}_{X_j}$ (i.e.~$f_i|_{X_j}$ agrees with an element of ${\cO}_{X_j}$ except at finitely many points of $X_j$).
\end{prop}

This, together with Corollary \ref{3.5} gives

\begin{thm}$K$ is $o$--minimal in ${\cL}$.
\end{thm}

\section{Further extensions}\label{extensions}
In this section we give extensions of the results of Sections 2 and
3 and the results of \cite{LR3}.

Let $G$ be an (additive) ordered abelian group.  Let $t$ be a symbol.  Then $t^G$ is a (multiplicative) ordered abelian group.  Following the notation of \cite{DMM1} and \cite{DMM2} (but not \cite{DMM3} or \cite{LR2}) we define ${\bR((t^G))}$ to be the maximally-complete valued field with additive value group $G$ (or multiplicative value group $t^G$) and residue field $\bR.$ So
$$
{\bR((t^G))} := \Big\{\sum_{g\in I} {a_{g}t^{g}}\colon a_g \in {\bR} \text{ and } I \subset G  \text{ well-ordered}\Big\}.
$$
We shall be a bit sloppy about mixing the additive and multiplicative valuations.  $I\subset G$ is well-ordered exactly when $t^I \subset t^G$ is reverse well-ordered.
The field $K$ of Puiseux series, or its completion, is a proper subfield of $\bR_1 := {\bR((t^\bQ))}.$
Considering $G = \bQ^m$ with the lexicographic ordering, we define
$$
\bR_m := {\bR((t^{\bQ^m}))}.
$$

It is clear that if $G_1 \subset G_2$ as ordered groups, then ${\bR((t^{G_1}))} \subset {\bR((t^{G_2}))}$ as valued fields.
Also, ${\bR((t^G))}$ is Henselian and, if $G$ is divisible, then ${\bR((t^G))}$ is real-closed.  We shall continue to use $<$ for the corresponding order on ${\bR((t^G))}.$

In analogy with Section 2, we define
\begin{defn}
\begin{eqnarray*}
{{\cal R}_{n,\alpha}(G)}&: =& A_{n,\alpha}{\widehat\otimes_{\bR}}^* {\bR((t^G))}\\
&:=&\Big\{\sum_{g\in I}{f_gt^g}\colon f_g \in A_{n,\alpha} \text{ and } I \subset G \text{ well ordered}\Big\}\\
{{\cal R}_n(G)}&:=&\bigcup_{\alpha>1}{{\cal R}_{n,\alpha}(G)}\\
{{\cal R}(G)}&: =&\bigcup_n\bigcup_{\alpha>1}{{\cal R}_{n,\alpha}(G)}.\\
\end{eqnarray*}
\end{defn}
As in Section 2, the elements of ${{\cal R}_n(G)}$ define functions from $I^n \subset {\bR((t^G))}^n$ to  ${\bR((t^G))}.$ Indeed
this interpretation is a ring endomorphism.  In other words, the field ${\bR((t^G))}$ has analytic ${{\cal R}(G)}$--structure.  (See \cite{CLR} and especially \cite{CL} for more about fields with analytic structure.)  The elements of ${{\cal R}_n(G)}$ are interpreted as zero on ${\bR((t^G))}^n \setminus I^n.$  Let
$$
{\cL}_G:=\langle +,\cdot,^{-1},0,1,<,{{\cal R}(G)}\rangle,
$$
so ${\bR((t^G))}$ is an ${\cL}_G$-structure.  Indeed, if $G_1 \subset G_2,$ then ${\bR((t^{G_2}))}$ is an ${\cL}_{G_1}$-structure.

\begin{thm}The Weierstrass Preparation Theorem (Theorem \ref{WPT}) and the Strong Noetherian Property (Theorem \ref{SNP}) hold with ${\cR}$ replaced by ${{\cal R}(G)}.$
\end{thm}
\begin{proof} Only minor modifications to the proofs of Theorems \ref{WPT} and \ref{SNP} are needed.
\end{proof}
The  arguments of Sections \ref{Notation&QE} and \ref{omin} show

\begin{thm}\label{GQE}
If $G$ is divisible then ${\bR((t^G))}$ admits quantifier
elimination and is $o$--minimal in ${\cL}_G.$
\end{thm}

\begin{cor}If $G_1 \subset G_2$ are divisible, then ${\bR((t^{G_1}))} \prec {\bR((t^{G_2}))}$ in ${\cL}_{G_1}.$
\end{cor}

We shall show in Section \ref{HP} that, though for $m<n$ we have
${\bR}_m \prec {\bR}_n$ in ${\cL}_{\bQ^m}$, there is a sentence of
${\cL}_{\bQ^n}$ that is true in ${\bR}_n$ but is not true in any
$o$--minimal expansion of ${\bR}_m .$

The results of \cite{LR3} also extend to this more general setting.

\begin{defn}We define the ring of \emph{*strictly convergent power series over ${\bR((t^G))}$} as
$$
{\bR((t^G))}^*\langle\xi\rangle := \{\sum_{g \in I}{a_g(\xi)t^g} \colon a_g(\xi) \in \bR[\xi] \text{ and $I$ well ordered}\},
$$
and the subring of \emph{*overconvergent power series over ${\bR((t^G))}$} as
\begin{eqnarray*}
{\bR((t^G))}^*\langle\langle\xi\rangle\rangle &:=& \{f  \colon f(\gamma\xi) \in {\bR((t^G))}^*\langle\xi\rangle \text{ for some } \gamma \in {\bR((t^G))}, \|\gamma\| >1\}, \\
{{\cal R}(G)}_{over} &:=& \bigcup_{n}{\bR((t^G))}^*\langle\langle\xi_1,\cdots, \xi_n\rangle\rangle,\\
\end{eqnarray*}
and the corresponding \emph{overconvergent language} as
$$
{\cL}_{G, over} := \langle +,\cdot,^{-1},0,1,<,{{\cal R}(G)}_{over}\rangle.
$$
\end{defn}

As in \cite{LR3} we have

\begin{thm}\label{GQEbis}
If G is divisible then ${\bR((t^G))}$ admits quantifier elimination and is $o$--minimal in ${\cL}_{G, over}$.
\end{thm}

Of course the $o$--minimality follows immediately from Theorem
\ref{GQE}.

\section{Extensions of the example of Hrushovski and Peterzil \cite{HP}}\label{HP}

In this section we show that with minor modifications, the idea of
the example of \cite{HP} can be iterated to give a nested family of
examples. This relates to a question of Hrushovski and Peterzil whether there
exists a small class of $o$--minimal structures such that any
sentence, true in some $o$--minimal structure, can be satisfied in
an expansion of a model in the class. Combining with expansions with
the exponential function, one perhaps can elaborate the tower of
examples further.

Consider the functional equation
\begin{equation}\tag{$\ast$}
F(\beta z)=\alpha z F(z)+1,
\end{equation}
and suppose that $F$ is a ``complex analytic'' solution for $|z|\leq 1$. By this we mean that, writing $z=x+\sqrt{-1}y$, $F(z)= f(x,y) + \sqrt{-1}g(x,y)$, $F(z)$ is differentiable as a function of $z$. This is a definable condition on the two ``real" functions, $f$, $g$ of the two ``real" variables $x$, $y$.
Then
$$
F(z)=\sum^\infty_{k=0} a_k z^k
$$
where
$$
a_k={\alpha^k\over \beta^{k(k+1)\over 2 }}.
$$
$(\alpha$, and $\beta$ are parameters).

By this we mean that for each $n \in \bN$ there is a constant $A_n$ such that
\begin{equation}\tag{$\ast\ast$}
|F(z)-\sum^n_{k=0} a_k z^k| \leq A_n |z^{n+1}|
\end{equation}
is true for all $z$ with $|z| \leq 1$.
Indeed, by \cite{PS} Theorem 2.50, one can take $A_n=C\cdot2^{n+1}$, for $C$ a constant independent of $n$.

Consider the following statement:\ $F(z)$ is a complex analytic function (in the above sense) on $|z|\leq 1$
 that satisfies $(*)$;
the number $\beta>0$ is within the radius of convergence of the
function $f(z)=\sum^\infty_{n=1}(n-1)!z^n$ and $\alpha >0$.

This statement is not satisfiable by any functions in any
$o$--minimal expansion of the field of Puiseux series $K_1$, or the
maximally complete field ${\bR}_1= {\bR((t^\bQ))},$ because, if it
were, we would have $\|\beta\|=\|t^\gamma\|,\
\|\alpha\|=\|t^\delta\|$, for some $0 < \gamma, \delta \in \bQ$, and
for suitable choice of $n$ the condition $(**)$ would be violated. On the other
hand, if we choose $\alpha,\beta\in {\bR}_2$ with $ord
(\alpha)=(1,0)$ and $ord(\beta)=(0,1)$ then $\Sigma a_k z^k\in
{\bR}_2\langle\langle z\rangle\rangle^*$ satisfies the statement on
${\bR}_2$.

This process can clearly be iterated to give, in the notation of Section \ref{extensions},
\begin{prop}For each $m$ there is a sentence of $\cL_{\bQ^m}$ true in ${\bR}_m$ but not satisfiable
in any $o$--minimal expansion of $\bR = {\bR}_0,
{\bR}_1,\ldots,{\bR}_{m-1}$.
\end{prop}

\end{document}